\documentclass[11pt]{article}
\usepackage{amssymb}

\def\hang{\hangindent\parindent}
\def\textindent#1{\indent\llap{#1\enspace}\ignorespaces}

\title{On the Calculation of gl.dim$G^{\mathbb{N}}(A)$ and gl.dim$\widetilde{A}$ by\\ Using Gr\"obner Bases
\thanks{Project supported by
the National Natural Science Foundation of China (10571038).\newline
e-mail: huishipp@yahoo.com}}

\author{Huishi Li\\
{\small Department of Applied Mathematics}\\
{\small College of Information Science and Technology}\\
{\small Hainan University}\\
{\small  Haikou 570228, China}}
\date{}

\begin{document}
\maketitle
\begin{center}
\begin{minipage}{120mm}
{\small {\bf Abstract.} Let $A=K\langle X_1,...,X_n\rangle /\langle
{\cal G}\rangle$ be a $K$-algebra defined by a finite Gr\"obner
basis ${\cal G}$. It is shown how to use the Ufnarovski graph
$\Gamma ({\bf LM}({\cal G}))$ and the graph of $n$-chains
$\Gamma_{\rm C}({\bf LM}({\cal G}))$ to calculate
gl.dim$G^{\mathbb{N}}(A)$ and gl.dim$\widetilde{A}$, where
$G^{\mathbb{N}}(A)$, respectively $\widetilde{A}$, is the associated
$\mathbb{N}$-graded algebra of $A$, respectively the Rees algebra of
$A$ with respect to the $\mathbb{N}$-filtration $FA$ of $A$ induced
by a weight $\mathbb{N}$-grading filtration of $K\langle
X_1,...,X_n\rangle$.  }
\end{minipage}\end{center}{\parindent=0pt\vskip 6pt
{\bf 2000 Mathematics Classification} Primary 16W70; Secondary 68W30
(16Z05).\vskip 6pt {\bf Key words} Global dimension, Growth of
algebra, Gr\"obner basis, $n$-chain}

%****************************************
\def\QED{\hfill{$\Box$}} \def\NZ{\mathbb{N}}
\def \r{\rightarrow}

\def\normalbaselines{\baselineskip 24pt\lineskip 4pt\lineskiplimit 4pt}
\def\mapdown#1{\llap{$\vcenter {\hbox {$\scriptstyle #1$}}$}
                                \Bigg\downarrow}
\def\mapdownr#1{\Bigg\downarrow\rlap{$\vcenter{\hbox
                                    {$\scriptstyle #1$}}$}}
\def\mapright#1#2{\smash{\mathop{\longrightarrow}\limits^{#1}_{#2}}}
\def\mapleft#1#2{\smash{\mathop{\longleftarrow}\limits^{#1}_{#2}}}
\def\mapup#1{\Bigg\uparrow\rlap{$\vcenter {\hbox  {$\scriptstyle #1$}}$}}
\def\mapupl#1{\llap{$\vcenter {\hbox {$\scriptstyle #1$}}$}
                                      \Bigg\uparrow}
\def\v5{\vskip .5truecm}
\def\T#1{\widetilde #1}
\def\OV#1{\overline {#1}}
\def\hang{\hangindent\parindent}
\def\textindent#1{\indent\llap{#1\enspace}\ignorespaces}
\def\item{\par\hang\textindent}
%****************************************
\message{<Paul Taylor's commutative diagrams, 20 July 1990>}
\newdimen\DiagramCellHeight\DiagramCellHeight3em %%
\newdimen\DiagramCellWidth\DiagramCellWidth3em %%
\newdimen\MapBreadth\MapBreadth.04em %%
\newdimen\MapShortFall\MapShortFall.4em %%
\newdimen\PileSpacing\PileSpacing1em %%
\def\labelstyle{\ifincommdiag\textstyle\else\scriptstyle\fi}%%
\let\objectstyle\displaystyle

%%======================================================================%

\def\rTo{\HorizontalMap\empty-\empty-\rhvee}%%
\def\lTo{\HorizontalMap\lhvee-\empty-\empty}%%
\def\dTo{\VerticalMap\empty|\empty|\dhvee}%%
\def\uTo{\VerticalMap\uhvee|\empty|\empty}%%
\let\uFrom\uTo\let\lFrom\lTo

\def\rArr{\HorizontalMap\empty-\empty-\rhla}%%
\def\lArr{\HorizontalMap\lhla-\empty-\empty}%%
\def\dArr{\VerticalMap\empty|\empty|\dhla}%%
\def\uArr{\VerticalMap\uhla|\empty|\empty}

\def\rDotsto{\HorizontalMap\empty\hfdot\hfdot\hfdot\rhvee}%%
\def\lDotsto{\HorizontalMap\lhvee\hfdot\hfdot\hfdot\empty}%%
\def\dDotsto{\VerticalMap\empty\vfdot\vfdot\vfdot\dhvee}%%
\def\uDotsto{\VerticalMap\uhvee\vfdot\vfdot\vfdot\empty}%%
\let\uDotsfrom\uDotsto\let\lDotsfrom\lDotsto

\def\rDashto{\HorizontalMap\empty\hfdash\hfdash\hfdash\rhvee}%%
\def\lDashto{\HorizontalMap\lhvee\hfdash\hfdash\hfdash\empty}%%
\def\dDashto{\VerticalMap\empty\vfdash\vfdash\vfdash\dhvee}%%
\def\uDashto{\VerticalMap\uhvee\vfdash\vfdash\vfdash\empty}%%
\let\uDashfrom\uDashto\let\lDashfrom\lDashto

\def\rImplies{\HorizontalMap==\empty=\Rightarrow}%%
\def\lImplies{\HorizontalMap\Leftarrow=\empty==}%%
\def\dImplies{\VerticalMap\|\|\empty\|\Downarrow}%%
\def\uImplies{\VerticalMap\Uparrow\|\empty\|\|}%%
\let\uImpliedby\uImplies\let\lImpliedby\lImplies

\def\rMapsto{\HorizontalMap\rtbar-\empty-\rhvee}%%
\def\lMapsto{\HorizontalMap\lhvee-\empty-\ltbar}%%
\def\dMapsto{\VerticalMap\dtbar|\empty|\dhvee}%%
\def\uMapsto{\VerticalMap\uhvee|\empty|\utbar}%%
\let\uMapsfrom\uMapsto\let\lMapsfrom\lMapsto

\def\rIntoA{\HorizontalMap\rthooka-\empty-\rhvee}%%
\def\rIntoB{\HorizontalMap\rthookb-\empty-\rhvee}%%
\def\lIntoA{\HorizontalMap\lhvee-\empty-\lthooka}%%
\def\lIntoB{\HorizontalMap\lhvee-\empty-\lthookb}%%
\def\dIntoA{\VerticalMap\dthooka|\empty|\dhvee}%%
\def\dIntoB{\VerticalMap\dthookb|\empty|\dhvee}%%
\def\uIntoA{\VerticalMap\uhvee|\empty|\uthooka}%%
\def\uIntoB{\VerticalMap\uhvee|\empty|\uthookb}%%
\let\uInfromA\uIntoA\let\uInfromB\uIntoB\let\lInfromA\lIntoA\let\lInfromB
\lIntoB\let\rInto\rIntoA\let\lInto\lIntoA\let\dInto\dIntoB\let\uInto\uIntoA

\def\rEmbed{\HorizontalMap\gt-\empty-\rhvee}%%
\def\lEmbed{\HorizontalMap\lhvee-\empty-\lt}%%
\def\dEmbed{\VerticalMap\vee|\empty|\dhvee}%%
\def\uEmbed{\VerticalMap\uhvee|\empty|\wedge}

\def\rProject{\HorizontalMap\empty-\empty-\triangleright}%%
\def\lProject{\HorizontalMap\triangleleft-\empty-\empty}%%
\def\uProject{\VerticalMap\triangleup|\empty|\empty}%%
\def\dProject{\VerticalMap\empty|\empty|\littletriangledown}

\def\rOnto{\HorizontalMap\empty-\empty-\twoheadrightarrow}%%
\def\lOnto{\HorizontalMap\twoheadleftarrow-\empty-\empty}%%
\def\dOnto{\VerticalMap\empty|\empty|\twoheaddownarrow}%%
\def\uOnto{\VerticalMap\twoheaduparrow|\empty|\empty}%%
\let\lOnfrom\lOnto\let\uOnfrom\uOnto

\def\hEq{\HorizontalMap==\empty==}%%
\def\vEq{\VerticalMap\|\|\empty\|\|}%%
\let\rEq\hEq\let\lEq\hEq\let\uEq\vEq\let\dEq\vEq

\def\hLine{\HorizontalMap\empty-\empty-\empty}%%
\def\vLine{\VerticalMap\empty|\empty|\empty}%%
\let\rLine\hLine\let\lLine\hLine\let\uLine\vLine\let\dLine\vLine

\def\hDots{\HorizontalMap\empty\hfdot\hfdot\hfdot\empty}%%
\def\vDots{\VerticalMap\empty\vfdot\vfdot\vfdot\empty}%%
\let\rDots\hDots\let\lDots\hDots\let\uDots\vDots\let\dDots\vDots

\def\hDashes{\HorizontalMap\empty\hfdash\hfdash\hfdash\empty}%%
\def\vDashes{\VerticalMap\empty\vfdash\vfdash\vfdash\empty}%%
\let\rDashes\hDashes\let\lDashes\hDashes\let\uDashes\vDashes\let\dDashes
\vDashes

\def\rPto{\HorizontalMap\empty-\empty-\rightharpoonup}%%
\def\lPto{\HorizontalMap\leftharpoonup-\empty-\empty}%%
\def\uPto{\VerticalMap\upharpoonright|\empty|\empty}%%
\def\dPto{\VerticalMap\empty|\empty|\downharpoonright}%%
\let\lPfrom\lPto\let\uPfrom\uPto

\def\NW{\NorthWest\DiagonalMap{\lah111}{\laf100}{}{\laf100}{}(2,2)}%%
\def\NE{\NorthEast\DiagonalMap{\lah22}{\laf0}{}{\laf0}{}(2,2)}%%
\def\SW{\SouthWest\DiagonalMap{}{\laf0}{}{\laf0}{\lah11}(2,2)}%%
\def\SE{\SouthEast\DiagonalMap{}{\laf100}{}{\laf100}{\lah122}(2,2)}

\def\nNW{\NorthWest\DiagonalMap{\lah135}{\laf112}{}{\laf112}{}(2,3)}%%
\def\nNE{\NorthEast\DiagonalMap{\lah36}{\laf12}{}{\laf12}{}(2,3)}%%
\def\sSW{\SouthWest\DiagonalMap{}{\laf12}{}{\laf12}{\lah35}(2,3)}%%
\def\sSE{\SouthEast\DiagonalMap{}{\laf112}{}{\laf112}{\lah136}(2,3)}

\def\wNW{\NorthWest\DiagonalMap{\lah153}{\laf121}{}{\laf121}{}(3,2)}%%
\def\eNE{\NorthEast\DiagonalMap{\lah63}{\laf21}{}{\laf21}{}(3,2)}%%
\def\wSW{\SouthWest\DiagonalMap{}{\laf21}{}{\laf21}{\lah53}(3,2)}%%
\def\eSE{\SouthEast\DiagonalMap{}{\laf121}{}{\laf121}{\lah163}(3,2)}

\def\NNW{\NorthWest\DiagonalMap{\lah113}{\laf101}{}{\laf101}{}(2,4)}%%
\def\NNE{\NorthEast\DiagonalMap{\lah25}{\laf01}{}{\laf01}{}(2,4)}%%
\def\SSW{\SouthWest\DiagonalMap{}{\laf01}{}{\laf01}{\lah13}(2,4)}%%
\def\SSE{\SouthEast\DiagonalMap{}{\laf101}{}{\laf101}{\lah125}(2,4)}

\def\WNW{\NorthWest\DiagonalMap{\lah131}{\laf110}{}{\laf110}{}(4,2)}%%
\def\ENE{\NorthEast\DiagonalMap{\lah52}{\laf10}{}{\laf10}{}(4,2)}%%
\def\WSW{\SouthWest\DiagonalMap{}{\laf10}{}{\laf10}{\lah31}(4,2)}%%
\def\ESE{\SouthEast\DiagonalMap{}{\laf110}{}{\laf110}{\lah152}(4,2)}

\def\NNNW{\NorthWest\DiagonalMap{\lah115}{\laf102}{}{\laf102}{}(2,6)}%%
\def\NNNE{\NorthEast\DiagonalMap{\lah16}{\laf02}{}{\laf02}{}(2,6)}%%
\def\SSSW{\SouthWest\DiagonalMap{}{\laf02}{}{\laf02}{\lah15}(2,6)}%%
\def\SSSE{\SouthEast\DiagonalMap{}{\laf102}{}{\laf102}{\lah116}(2,6)}

\def\WWNW{\NorthWest\DiagonalMap{\lah151}{\laf120}{}{\laf120}{}(6,2)}%%
\def\EENE{\NorthEast\DiagonalMap{\lah61}{\laf20}{}{\laf20}{}(6,2)}%%
\def\WWSW{\SouthWest\DiagonalMap{}{\laf20}{}{\laf20}{\lah51}(6,2)}%%
\def\EESE{\SouthEast\DiagonalMap{}{\laf120}{}{\laf120}{\lah161}(6,2)}

%%======================================================================%

\font\tenln=line10

\mathchardef\lt="313C \mathchardef\gt="313E

\def\rhvee{\mkern-10mu\gt}%%
\def\lhvee{\lt\mkern-10mu}%%
\def\dhvee{\vbox\tozpt{\vss\hbox{$\vee$}\kern0pt}}%%
\def\uhvee{\vbox\tozpt{\hbox{$\wedge$}\vss}}%%
\def\rhcvee{\mkern-10mu\succ}%%
\def\lhcvee{\prec\mkern-10mu}%%
\def\dhcvee{\vbox\tozpt{\vss\hbox{$\curlyvee$}\kern0pt}}%%
\def\uhcvee{\vbox\tozpt{\hbox{$\curlywedge$}\vss}}%%
\def\rhvvee{\mkern-10mu\gg}%%
\def\lhvvee{\ll\mkern-10mu}%%
\def\dhvvee{\vbox\tozpt{\vss\hbox{$\vee$}\kern-.6ex\hbox{$\vee$}\kern0pt}}%%
\def\uhvvee{\vbox\tozpt{\hbox{$\wedge$}\kern-.6ex\hbox{$\wedge$}\vss}}%%
\def\twoheaddownarrow{\rlap{$\downarrow$}\raise-.5ex\hbox{$\downarrow$}}%%
\def\twoheaduparrow{\rlap{$\uparrow$}\raise.5ex\hbox{$\uparrow$}}%%
\def\triangleup{{\scriptscriptstyle\bigtriangleup}}%%
\def\littletriangledown{{\scriptscriptstyle\triangledown}}%%
\def\rhla{\vbox\tozpt{\vss\hbox\tozpt{\hss\tenln\char'55}\kern\axisheight}}%%
\def\lhla{\vbox\tozpt{\vss\hbox\tozpt{\tenln\char'33\hss}\kern\axisheight}}%%
\def\dhla{\vbox\tozpt{\vss\hbox\tozpt{\tenln\char'77\hss}}}%%
\def\uhla{\vbox\tozpt{\hbox\tozpt{\tenln\char'66\hss}\vss}}%%
\def\htdot{\mkern3.15mu\cdot\mkern3.15mu}%%
\def\vtdot{\vbox to 1.46ex{\vss\hbox{$\cdot$}}}%%
\def\utbar{\vrule height 0.093ex depth0pt width 0.4em} \let\dtbar\utbar%%
\def\rtbar{\mkern1.5mu\vrule height 1.1ex depth.06ex width .04em\mkern1.5mu}%
\let\ltbar\rtbar%%
\def\rthooka{\raise.603ex\hbox{$\scriptscriptstyle\subset$}}%%
\def\lthooka{\raise.603ex\hbox{$\scriptscriptstyle\supset$}}%%
\def\rthookb{\raise-.022ex\hbox{$\scriptscriptstyle\subset$}}%%
\def\lthookb{\raise-.022ex\hbox{$\scriptscriptstyle\supset$}}%%
\def\dthookb{\hbox{$\scriptscriptstyle\cap$}\mkern5.5mu}%%
\def\uthookb{\hbox{$\scriptscriptstyle\cup$}\mkern4.5mu}%%
\def\dthooka{\mkern6mu\hbox{$\scriptscriptstyle\cap$}}%%
\def\uthooka{\mkern6mu\hbox{$\scriptscriptstyle\cup$}}%%
\def\hfdot{\mkern3.15mu\cdot\mkern3.15mu}%%
\def\vfdot{\vbox to 1.46ex{\vss\hbox{$\cdot$}}}%%
\def\vfdashstrut{\vrule width0pt height1.3ex depth0.7ex}%%
\def\vfthedash{\vrule width\MapBreadth height0.6ex depth 0pt}%%
\def\hfthedash{\vrule\horizhtdp width 0.26em}%%
\def\hfdash{\mkern5.5mu\hfthedash\mkern5.5mu}%%
\def\vfdash{\vfdashstrut\vfthedash}%%

\def\nwhTO{\nwarrow\mkern-1mu}%%
\def\nehTO{\mkern-.1mu\nearrow}%%
\def\sehTO{\searrow\mkern-.02mu}%%
\def\swhTO{\mkern-.8mu\swarrow}%%

\def\SEpbk{\rlap{\smash{\kern0.1em \vrule depth 2.67ex height -2.55ex width 0%
.9em \vrule height -0.46ex depth 2.67ex width .05em }}}%%
\def\SWpbk{\llap{\smash{\vrule height -0.46ex depth 2.67ex width .05em \vrule
depth 2.67ex height -2.55ex width .9em \kern0.1em }}}%%
\def\NEpbk{\rlap{\smash{\kern0.1em \vrule depth -3.48ex height 3.67ex width 0%
.95em \vrule height 3.67ex depth -1.39ex width .05em }}}%%
\def\NWpbk{\llap{\smash{\vrule height 3.6ex depth -1.39ex width .05em \vrule
depth -3.48ex height 3.67ex width .95em \kern0.1em }}}

%%======================================================================%

\newcount\cdna\newcount\cdnb\newcount\cdnc\newcount\cdnd\cdna=\catcode`\@%
\catcode`\@=11 \let\then\relax\def\loopa#1\repeat{\def\bodya{#1}\iteratea}%
\def\iteratea{\bodya\let\next\iteratea\else\let\next\relax\fi\next}\def\loopb
#1\repeat{\def\bodyb{#1}\iterateb}\def\iterateb{\bodyb\let\next\iterateb\else
\let\next\relax\fi\next} \def\mapctxterr{\message{commutative diagram: map
context error}}\def\mapclasherr{\message{commutative diagram: clashing maps}}%
\def\ObsDim#1{\expandafter\message{! diagrams Warning: Dimension \string#1 is
obsolete
(ignored)}\global\let#1\ObsDimq\ObsDimq}\def\ObsDimq{\dimen@=}\def
\HorizontalMapLength{\ObsDim\HorizontalMapLength}\def\VerticalMapHeight{%
\ObsDim\VerticalMapHeight}\def\VerticalMapDepth{\ObsDim\VerticalMapDepth}\def
\VerticalMapExtraHeight{\ObsDim\VerticalMapExtraHeight}\def
\VerticalMapExtraDepth{\ObsDim\VerticalMapExtraDepth}\def\ObsCount#1{%
\expandafter\message{! diagrams Warning: Count \string#1 is obsolete (ignored%
)}\global\let#1\ObsCountq\ObsCountq}\def\ObsCountq{\count@=}\def
\DiagonalLineSegments{\ObsCount\DiagonalLineSegments}\def\tozpt{to\z@}\def
\sethorizhtdp{\dimen8=\axisheight\dimen9=\MapBreadth\advance\dimen8.5\dimen9%
\advance\dimen9-\dimen8}\def\horizhtdp{height\dimen8 depth\dimen9
}\def \axisheight{\fontdimen22\the\textfont2 }\countdef\boxc@unt=14

\def\bombparameters{\hsize\z@\rightskip\z@ plus1fil minus\maxdimen
\parfillskip\z@\linepenalty9000 \looseness0 \hfuzz\maxdimen\hbadness10000
\clubpenalty0 \widowpenalty0 \displaywidowpenalty0
\interlinepenalty0 \predisplaypenalty0 \postdisplaypenalty0
\interdisplaylinepenalty0
\interfootnotelinepenalty0 \floatingpenalty0 \brokenpenalty0 \everypar{}%
\leftskip\z@\parskip\z@\parindent\z@\pretolerance10000
\tolerance10000 \hyphenpenalty10000 \exhyphenpenalty10000
\binoppenalty10000 \relpenalty10000 \adjdemerits0
\doublehyphendemerits0 \finalhyphendemerits0 \prevdepth\z@}\def
\startbombverticallist{\hbox{}\penalty1\nointerlineskip}

\def\pushh#1\to#2{\setbox#2=\hbox{\box#1\unhbox#2}}\def\pusht#1\to#2{\setbox#%
2=\hbox{\unhbox#2\box#1}}

\newif\ifallowhorizmap\allowhorizmaptrue\newif\ifallowvertmap
\allowvertmapfalse\newif\ifincommdiag\incommdiagfalse

\def\diagram{\hbox\bgroup$\vcenter\bgroup\startbombverticallist
\incommdiagtrue\baselineskip\DiagramCellHeight\lineskip\z@\lineskiplimit\z@
\mathsurround\z@\tabskip\z@\let\\\diagcr\allowhorizmaptrue\allowvertmaptrue
\halign\bgroup\lcdtempl##\rcdtempl&&\lcdtempl##\rcdtempl\cr}\def\enddiagram{%
\crcr\egroup\reformatmatrix\egroup$\egroup}\def\commdiag#1{{\diagram#1%
\enddiagram}}

\def\lcdtempl{\futurelet\thefirsttoken\dolcdtempl}\newif\ifemptycell\def
\dolcdtempl{\ifx\thefirsttoken\rcdtempl\then\hskip1sp plus 1fil
\emptycelltrue
\else\hfil$\emptycellfalse\objectstyle\fi}\def\rcdtempl{\ifemptycell\else$%
\hfil\fi}\def\diagcr{\cr} \def\across#1{\span\omit\mscount=#1
\loop\ifnum \mscount>2
\spAn\repeat\ignorespaces}\def\spAn{\relax\span\omit\advance
\mscount by -1}

\def\CellSize{\afterassignment\cdhttowd\DiagramCellHeight}\def\cdhttowd{%
\DiagramCellWidth\DiagramCellHeight}\def\MapsAbut{\MapShortFall\z@}

\newcount\cdvdl\newcount\cdvdr\newcount\cdvd\newcount\cdbfb\newcount\cdbfr
\newcount\cdbfl\newcount\cdvdr\newcount\cdvdl\newcount\cdvd

\def\reformatmatrix{\bombparameters\cdvdl=\insc@unt\cdvdr=\cdvdl\cdbfb=%
\boxc@unt\advance\cdbfb1
\cdbfr=\cdbfb\setbox1=\vbox{}\dimen2=\z@\loop\setbox
0=\lastbox\ifhbox0 \dimen1=\lastskip\unskip\dimen5=\ht0
\advance\dimen5 \dimen
1 \dimen4=\dp0 \penalty1 \reformatrow\unpenalty\ht4=\dimen5 \dp4=\dimen4 \ht3%
\z@\dp3\z@\setbox1=\vbox{\box4 \nointerlineskip\box3 \nointerlineskip\unvbox1%
}\dimen2=\dimen1 \repeat\unvbox1}

\newif\ifcontinuerow

\def\reformatrow{\cdbfl=\cdbfr\noindent\unhbox0 \loopa\unskip\setbox\cdbfl=%
\lastbox\ifhbox\cdbfl\advance\cdbfl1\repeat\par\unskip\dimen6=2%
\DiagramCellWidth\dimen7=-\DiagramCellWidth\setbox3=\hbox{}\setbox4=\hbox{}%
\setbox7=\box\voidb@x\cdvd=\cdvdl\continuerowtrue\loopa\advance\cdvd-1
\adjustcells\ifcontinuerow\advance\dimen6\wd\cdbfl\cdda=.5\dimen6
\ifdim\cdda
<\DiagramCellWidth\then\dimen6\DiagramCellWidth\advance\dimen6-\cdda
\nopendvert\cdda\DiagramCellWidth\fi\advance\dimen7\cdda\dimen6=\wd\cdbfl
\reformatcell\advance\cdbfl-1 \repeat\advance\dimen7.5\dimen6
\outHarrow} \def
\adjustcells{\ifnum\cdbfr>\cdbfl\then\ifnum\cdvdr>\cdvd\then\continuerowfalse
\else\setbox\cdbfl=\hbox
to\wd\cdvd{\lcdtempl\VonH{}\rcdtempl}\fi\else\ifnum
\cdvdr>\cdvd\then\advance\cdvdr-1
\setbox\cdvd=\vbox{}\wd\cdvd=\wd\cdbfl\dp \cdvd=\dp1 \fi\fi}

\def\reformatcell{\sethorizhtdp\noindent\unhbox\cdbfl\skip0=\lastskip\unskip
\par\ifcase\prevgraf\reformatempty\or\reformatobject\else\reformatcomplex\fi
\unskip}\def\reformatobject{\setbox6=\lastbox\unskip\vadjdon6\outVarrow
\setbox6=\hbox{\unhbox6}\advance\dimen7-.5\wd6
\outHarrow\dimen7=-.5\wd6 \pusht6\to4}\newcount\globnum

\def\reformatcomplex{\setbox6=\lastbox\unskip\setbox9=\lastbox\unskip\setbox9%
=\hbox{\unhbox9
\skip0=\lastskip\unskip\global\globnum\lastpenalty\hskip\skip 0
}\advance\globnum9999
\ifcase\globnum\reformathoriz\or\reformatpile\or
\reformatHonV\or\reformatVonH\or\reformatvert\or\reformatHmeetV\fi}

\def\reformatempty{\vpassdon\ifdim\skip0>\z@\then\hpassdon\else\ifvoid2 \then
\else\advance\dimen7-.5\dimen0
\cdda=\wd2\advance\cdda.5\dimen0\wd2=\cdda\fi
\fi}\def\VonH{\doVonH6}\def\HonV{\doVonH7}\def\HmeetV{\MapBreadth-2%
\MapShortFall\doVonH4}\def\doVonH#1{\cdna-999#1\futurelet\thenexttoken
\dooVonH}\def\dooVonH{\let\next\relax\sethorizhtdp\ifallowhorizmap
\ifallowvertmap\then\ifx\thenexttoken[\then\let\next\VonHstrut\else
\sethorizhtdp\dimen0\MapBreadth\let\next\VonHnostrut\fi\else\mapctxterr\fi
\else\mapctxterr\fi\next}\def\VonHstrut[#1]{\setbox0=\hbox{$#1$}\dimen0\wd0%
\dimen8\ht0\dimen9\dp0 \VonHnostrut}\def\VonHnostrut{\setbox0=\hbox{}\ht0=%
\dimen8\dp0=\dimen9\wd0=.5\dimen0 \copy0\penalty\cdna\box0
\allowhorizmapfalse
\allowvertmapfalse}\def\reformatHonV{\hpassdon\doreformatHonV}\def
\reformatHmeetV{\dimen@=\wd9 \advance\dimen7-\wd9 \outHarrow\setbox6=\hbox{%
\unhbox6}\dimen7-\wd6 \advance\dimen@\wd6 \setbox6=\hbox to\dimen@{\hss}%
\pusht6\to4\doreformatHonV}\def\doreformatHonV{\setbox9=\hbox{\unhbox9
\unskip
\unpenalty\global\setbox\globbox=\lastbox}\vadjdon\globbox\outVarrow}\def
\reformatVonH{\vpassdon\advance\dimen7-\wd9 \outHarrow\setbox6=\hbox{\unhbox6%
}\dimen7=-\wd6 \setbox6=\hbox{\kern\wd9 \kern\wd6}\pusht6\to4}\def\hpassdon{}%
\def\vpassdon{\dimen@=\dp\cdvd\advance\dimen@\dimen4 \advance\dimen@\dimen5
\dp\cdvd=\dimen@\nopendvert}\def\vadjdon#1{\dimen8=\ht#1
\dimen9=\dp#1 }

\def\HorizontalMap#1#2#3#4#5{\sethorizhtdp\setbox1=\makeharrowpart{#1}\def
\arrowfillera{#2}\def\arrowfillerb{#4}\setbox5=\makeharrowpart{#5}\ifx
\arrowfillera\justhorizline\then\def\arra{\hrule\horizhtdp}\def\kea{\kern-0.%
01em}\let\arrstruthtdp\horizhtdp\else\def\kea{\kern-0.15em}\setbox2=\hbox{%
\kea${\arrowfillera}$\kea}\def\arra{\copy2}\def\arrstruthtdp{height\ht2 depth%
\dp2
}\fi\ifx\arrowfillerb\justhorizline\then\def\arrb{\hrule\horizhtdp}\def
\keb{kern-0.01em}\ifx\arrowfillera\empty\then\let\arrstruthtdp\horizhtdp\fi
\else\def\keb{\kern-0.15em}\setbox4=\hbox{\keb${\arrowfillerb}$\keb}\def\arrb
{\copy4}\ifx\arrowfilera\empty\then\def\arrstruthtdp{height\ht4 depth\dp4 }%
\fi\fi\setbox3=\makeharrowpart{{#3}\vrule width\z@\arrstruthtdp}%
\ifallowhorizmap\then\let\execmap\execHorizontalMap\else\let\execmap
\mapctxterr\fi\allowhorizmapfalse\gettwoargs}\def\makeharrowpart#1{\hbox{%
\mathsurround\z@\edef\next{#1}\ifx\next\empty\else$\mkern-1.5mu{\next}\mkern-%
1.5mu$\fi}}\def\justhorizline{-}

\def\execHorizontalMap{\dimen0=\wd6 \ifdim\dimen0<\wd7\then\dimen0=\wd7\fi
\dimen3=\wd3 \ifdim\dimen0<2em\then\dimen0=2em\fi\skip2=.5\dimen0
\ifincommdiag plus 1fill\fi minus\z@\advance\skip2-.5\dimen3
\skip4=\skip2 \advance\skip2-\wd1 \advance\skip4-\wd5
\kern\MapShortFall\box1 \xleaders
\arra\hskip\skip2 \vbox{\lineskiplimit\maxdimen\lineskip.5ex \ifhbox6 \hbox to%
\dimen3 {\hss\box6\hss}\fi\vtop{\box3 \ifhbox7 \hbox to\dimen3
{\hss\box7\hss
}\fi}}\ifincommdiag\kern-.5\dimen3\penalty-9999\null\kern.5\dimen3\fi
\xleaders\arrb\hskip\skip4 \box5 \kern\MapShortFall}

\def\reformathoriz{\vadjdon6\outVarrow\ifvoid7\else\mapclasherr\fi\setbox2=%
\box9 \wd2=\dimen7 \dimen7=\z@\setbox7=\box6 }

\def\resetharrowpart#1#2{\ifvoid#1\then\ifdim#2=\z@\else\setbox4=\hbox{%
\unhbox4\kern#2}\fi\else\ifhbox#1\then\setbox#1=\hbox
to#2{\unhbox#1}\else
\widenpile#1\fi\pusht#1\to4\fi}\def\outHarrow{\resetharrowpart2{\wd2}\pusht2%
\to4\resetharrowpart7{\dimen7}\pusht7\to4\dimen7=\z@}

\def\pile#1{{\incommdiagtrue\let\pile\innerpile\allowvertmapfalse
\allowhorizmaptrue\baselineskip.5\PileSpacing\lineskip\z@\lineskiplimit\z@
\mathsurround\z@\tabskip\z@\let\\\pilecr\vcenter{\halign{\hfil$##$\hfil\cr#1
\crcr}}}\ifincommdiag\then\ifallowhorizmap\then\penalty-9998
\allowvertmapfalse\allowhorizmapfalse\else\mapctxterr\fi\fi}\def\pilecr{\cr}%
\def\innerpile#1{\noalign{\halign{\hfil$##$\hfil\cr#1 \crcr}}}

\def\reformatpile{\vadjdon9\outVarrow\ifvoid7\else\mapclasherr\fi\penalty1
\setbox9=\hbox{\unhbox9 \unskip\unpenalty\setbox9=\lastbox\unhbox9
\global
\setbox\globbox=\lastbox}\unvbox\globbox\setbox9=\vbox{}\setbox7=\vbox{}%
\loopb\setbox6=\lastbox\ifhbox6
\skip3=\lastskip\unskip\splitpilerow\repeat
\unpenalty\setbox9=\hbox{$\vcenter{\unvbox9}$}\setbox2=\box9
\dimen7=\z@}\def
\pilestrut{\vrule height\dimen0 depth\dimen3 width\z@}\def\splitpilerow{%
\dimen0=\ht6 \dimen3=\dp6
\noindent\unhbox6\unskip\setbox6=\lastbox\unskip
\unhbox6\par\setbox6=\lastbox\unskip\ifcase\prevgraf\or\setbox6=\hbox\tozpt{%
\hss\unhbox6\hss}\ht6=\dimen0 \dp6=\dimen3
\setbox9=\vbox{\vskip\skip3 \hbox
to\dimen7{\hfil\box6}\nointerlineskip\unvbox9}\setbox7=\vbox{\vskip\skip3
\hbox{\pilestrut\hfil}\nointerlineskip\unvbox7}\or\setbox7=\vbox{\vskip\skip3
\hbox{\pilestrut\unhbox6}\nointerlineskip\unvbox7}\setbox6=\lastbox\unskip
\setbox9=\vbox{\vskip\skip3 \hbox to\dimen7{\pilestrut\unhbox6}%
\nointerlineskip\unvbox9}\fi\unskip}

\def\widenpile#1{\setbox#1=\hbox{$\vcenter{\unvbox#1 \setbox8=\vbox{}\loopb
\setbox9=\lastbox\ifhbox9
\skip3=\lastskip\unskip\setbox8=\vbox{\vskip\skip3 \hbox
to\dimen7{\unhbox9}\nointerlineskip\unvbox8}\repeat\unvbox8 }$}}

\def\justverticalline{|}\def\makevarrowpart#1{\hbox to\MapBreadth{\hss$\kern
\MapBreadth{#1}$\hss}}\def\VerticalMap#1#2#3#4#5{\setbox1=\makevarrowpart{#1}%
\def\arrowfillera{#2}\setbox3=\makevarrowpart{#3}\def\arrowfillerb{#4}\setbox
5=\makevarrowpart{#5}\ifx\arrowfillera\justverticalline\then\def\arra{\vrule
width\MapBreadth}\def\kea{\kern-0.05ex}\else\def\kea{\kern-0.35ex}\setbox2=%
\vbox{\kea\makevarrowpart\arrowfillera\kea}\def\arra{\copy2}\fi\ifx
\arrowfillerb\justverticalline\then\def\arrb{\vrule
width\MapBreadth}\def\keb
{\kern-0.05ex}\else\def\keb{\kern-0.35ex}\setbox4=\vbox{\keb\makevarrowpart
\arrowfillerb\keb}\def\arrb{\copy4}\fi\ifallowvertmap\then\let\execmap
\execVerticalMap\else\let\execmap\mapctxterr\fi\allowhorizmapfalse\gettwoargs
}

\def\execVerticalMap{\setbox3=\makevarrowpart{\box3}\setbox0=\hbox{}\ht0=\ht3
\dp0\z@\ht3\z@\box6 \setbox8=\vtop spread2ex{\offinterlineskip\box3
\xleaders
\arrb\vfill\box5 \kern\MapShortFall}\dp8=\z@\box8 \kern-\MapBreadth\setbox8=%
\vbox spread2ex{\offinterlineskip\kern\MapShortFall\box1
\xleaders\arra\vfill \box0}\ht8=\z@\box8
\ifincommdiag\then\kern-.5\MapBreadth\penalty-9995 \null
\kern.5\MapBreadth\fi\box7\hfil}

\newcount\colno\newdimen\cdda\newbox\globbox\def\reformatvert{\setbox6=\hbox{%
\unhbox6}\cdda=\wd6 \dimen3=\dp\cdvd\advance\dimen3\dimen4
\setbox\cdvd=\hbox {}\colno=\prevgraf\advance\colno-2
\loopb\setbox9=\hbox{\unhbox9 \unskip
\unpenalty\dimen7=\lastkern\unkern\global\setbox\globbox=\lastbox\advance
\dimen7\wd\globbox\advance\dimen7\lastkern\unkern\setbox9=\lastbox\vtop to%
\dimen3{\unvbox9}\kern\dimen7 }\ifnum\colno>0
\ifdim\wd9<\PileSpacing\then \setbox9=\hbox
to\PileSpacing{\unhbox9}\fi\dimen0=\wd9 \advance\dimen0-\wd
\globbox\setbox\cdvd=\hbox{\kern\dimen0 \box\globbox\unhbox\cdvd}\pushh9\to6%
\advance\colno-1
\setbox9=\lastbox\unskip\repeat\advance\dimen7-.5\wd6
\advance\dimen7.5\cdda\advance\dimen7-\wd9 \outHarrow\dimen7=-.5\wd6
\advance
\dimen7-.5\cdda\pusht9\to4\pusht6\to4\nopendvert\dimen@=\dimen6\advance
\dimen@-\wd\cdvd\advance\dimen@-\wd\globbox\divide\dimen@2
\setbox\cdvd=\hbox
{\kern\dimen@\box\globbox\unhbox\cdvd\kern\dimen@}\dimen8=\dp\cdvd\advance
\dimen8\dimen5 \dp\cdvd=\dimen8 \ht\cdvd=\z@}

\def\outVarrow{\ifhbox\cdvd\then\deepenbox\cdvd\pusht\cdvd\to3\else
\nopendvert\fi\dimen3=\dimen5 \advance\dimen3-\dimen8 \setbox\cdvd=\vbox{%
\vfil}\dp\cdvd=\dimen3} \def\nopendvert{\setbox3=\hbox{\unhbox3\kern\dimen6}}%
\def\deepenbox\cdvd{\setbox\cdvd=\hbox{\dimen3=\dimen4 \advance\dimen3-\dimen
9 \setbox6=\hbox{}\ht6=\dimen3 \dp6=-\dimen3 \dimen0=\dp\cdvd\advance\dimen0%
\dimen3 \unhbox\cdvd\dimen3=\lastkern\unkern\setbox8=\hbox{\kern\dimen3}%
\loopb\setbox9=\lastbox\ifvbox9 \setbox9=\vtop to\dimen0{\copy6
\nointerlineskip\unvbox9 }\dimen3=\lastkern\unkern\setbox8=\hbox{\kern\dimen3%
\box9\unhbox8}\repeat\unhbox8 }}

\newif\ifPositiveGradient\PositiveGradienttrue\newif\ifClimbing\Climbingtrue
\newcount\DiagonalChoice\DiagonalChoice1 \newcount\lineno\newcount\rowno
\newcount\charno\def\laf{\afterassignment\xlaf\charno='}\def\xlaf{\hbox{%
\tenln\char\charno}}\def\lah{\afterassignment\xlah\charno='}\def\xlah{\hbox{%
\tenln\char\charno}}\def\makedarrowpart#1{\hbox{\mathsurround\z@${#1}$}}\def
\lad{\afterassignment\xlad\charno='}\def\xlad{\setbox2=\xlaf\setbox0=\hbox to%
.5\wd2{$\hss\ldot\hss$}\ht0=.25\ht2 \dp0=\ht0 \hbox{\mv-\ht0\copy0 \mv\ht0%
\box0}}

\def\DiagonalMap#1#2#3#4#5{\ifPositiveGradient\then\let\mv\raise\else\let\mv
\lower\fi\setbox2=\makedarrowpart{#2}\setbox1=\makedarrowpart{#1}\setbox4=%
\makedarrowpart{#4}\setbox5=\makedarrowpart{#5}\setbox3=\makedarrowpart{#3}%
\let\execmap\execDiagonalLine\gettwoargs}

\def\makeline#1(#2,#3;#4){\hbox{\dimen1=#2\relax\dimen2=#3\relax\dimen5=#4%
\relax\vrule height\dimen5 depth\z@ width\z@\setbox8=\hbox to\dimen1{\tenln#1%
\hss}\cdna=\dimen5 \divide\cdna\dimen2 \ifnum\cdna=0 \then\box8 \else\dimen4=%
\dimen5 \advance\dimen4-\dimen2 \divide\dimen4\cdna\dimen3=\dimen1 \cdnb=%
\dimen2 \divide\cdnb1000 \divide\dimen3\cdnb\cdnb=\dimen4
\divide\cdnb1000 \multiply\dimen3\cdnb\dimen6\dimen1
\advance\dimen6-\dimen3 \cdnb=0
\ifPositiveGradient\then\dimen7\z@\else\dimen7\cdna\dimen4
\multiply\dimen4-1 \fi\loop\raise\dimen7\copy8
\ifnum\cdnb<\cdna\hskip-\dimen6 \advance\cdnb1
\advance\dimen7\dimen4
\repeat\fi}}\newdimen\objectheight\objectheight1.5ex

\def\execDiagonalLine{\setbox0=\hbox\tozpt{\cdna=\xcoord\cdnb=\ycoord\dimen8=%
\wd2 \dimen9=\ht2 \dimen0=\cdnb\DiagramCellHeight\advance\dimen0-2%
\MapShortFall\advance\dimen0-\objectheight\setbox2=\makeline\box2(\dimen8,%
\dimen9;.5\dimen0)\setbox4=\makeline\box4(\dimen8,\dimen9;.5\dimen0)\dimen0=2%
\wd2 \advance\dimen0-\cdna\DiagramCellWidth\advance\dimen0
2\DiagramCellWidth
\dimen2\DiagramCellHeight\advance\dimen2-\MapShortFall\dimen1\dimen2
\advance \dimen1-\ht1 \advance\dimen2-\ht2 \dimen6=\dimen2
\advance\dimen6.25\dimen8 \dimen3\dimen2 \advance\dimen3-\ht3
\dimen4=\dimen2 \dimen7=\dimen2 \advance \dimen4-\ht4
\advance\dimen7-\ht7 \advance\dimen7-.25\dimen8
\ifPositiveGradient\then\hss\raise\dimen4\hbox{\rlap{\box5}\box4}\llap{\raise
\dimen6\box6\kern.25\dimen9}\else\kern-.5\dimen0 \rlap{\raise\dimen1\box1}%
\raise\dimen2\box2 \llap{\raise\dimen7\box7\kern.25\dimen9}\fi\raise\dimen3%
\hbox\tozpt{\hss\box3\hss}\ifPositiveGradient\then\rlap{\kern.25\dimen9\raise
\dimen7\box7}\raise\dimen2\box2\llap{\raise\dimen1\box1}\kern-.5\dimen0
\else
\rlap{\kern.25\dimen9\raise\dimen6\box6}\raise\dimen4\hbox{\box4\llap{\box5}}%
\hss\fi}\ht0\z@\dp0\z@\box0}

\def\NorthWest{\PositiveGradientfalse\Climbingtrue\DiagonalChoice0 }\def
\NorthEast{\PositiveGradienttrue\Climbingtrue\DiagonalChoice1
}\def\SouthWest
{\PositiveGradienttrue\Climbingfalse\DiagonalChoice3 }\def\SouthEast{%
\PositiveGradientfalse\Climbingfalse\DiagonalChoice2 }

\newif\ifmoremapargs\def\gettwoargs{\setbox7=\box\voidb@x\setbox6=\box
\voidb@x\moremapargstrue\def\whichlabel{6}\def\xcoord{2}\def\ycoord{2}\def
\contgetarg{\def\whichlabel{7}\ifmoremapargs\then\let\next\getanarg\let
\contgetarg\execmap\else\let\next\execmap\fi\next}\getanarg}\def\getanarg{%
\futurelet\thenexttoken\switcharg}\def\getlabel#1#2#3{\setbox#1=\hbox{$%
\labelstyle\>{#3}\>$}\dimen0=\ht#1\advance\dimen0 .4ex\ht#1=\dimen0 \dimen0=%
\dp#1\advance\dimen0 .4ex\dp#1=\dimen0 \contgetarg}\def\eatspacerepeat{%
\afterassignment\getanarg\let\junk=
}\def\catcase#1:{{\ifcat\noexpand
\thenexttoken#1\then\global\let\xcase\docase\fi}\xcase}\def\tokcase#1:{{\ifx
\thenexttoken#1\then\global\let\xcase\docase\fi}\xcase}\def\default:{\docase}%
\def\docase#1\esac#2\esacs{#1}\def\skipcase#1\esac{}\def\getcoordsrepeat(#1,#%
2){\def\xcoord{#1}\def\ycoord{#2}\getanarg}\let\esacs\relax\def\switcharg{%
\global\let\xcase\skipcase\catcase{&}:\moremapargsfalse\contgetarg\esac
\catcase\bgroup:\getlabel\whichlabel-\esac\catcase^:\getlabel6\esac\catcase_:%
\getlabel7\esac\tokcase{~}:\getlabel3\esac\tokcase(:\getcoordsrepeat\esac
\catcase{
}:\eatspacerepeat\esac\default:\moremapargsfalse\contgetarg\esac
\esacs}

\catcode`\@=\cdna
%****************************************

\def\LH{{\bf LH}}\def\LM{{\bf LM}}\def\LT{{\bf
LT}}\def\KS{K\langle X\rangle}
\def\B{{\cal B}} \def\LC{{\bf LC}} \def\G{{\cal G}} \def\FRAC#1#2{\displaystyle{\frac{#1}{#2}}}
\def\SUM^#1_#2{\displaystyle{\sum^{#1}_{#2}}} \def\O{{\cal O}}  \def\J{{\bf J}}
{\parindent=0pt\vskip 1truecm Let $\KS =K\langle X_1,...,X_n\rangle$
be the free algebra generated by $X=\{ X_1,...,X_n\}$ over a field
$K$, $I$ an arbitrary (two-sided) ideal of $\KS$, and $A=\KS /I$ the
corresponding quotient algebra. Fixing  a positive weight
$\NZ$-gradation $\{ \KS_p\}_{p\in\NZ}$ for $\KS$ by assigning to
each $X_i$ a positive degree $n_i$, $1\le i\le n$,  so that $\KS
=\oplus_{p\in\NZ}\KS_p$, and considering the $\NZ$-filtration $FA$
of $A$ induced by the weight $\NZ$-grading filtration $F\KS=\{
F_p\KS=\oplus_{i\le p}\KS_i\}_{p\in\NZ}$ of $\KS$, then it is
well-known that $FA$ determines two $\NZ$-graded $K$-algebras,
namely  the associated $\NZ$-graded algebra
$G^{\NZ}(A)=\oplus_{p\in\NZ}(F_pA/F_{p-1}A)$ of $A$ and the Rees
algebra $\T A=\oplus_{p\in\NZ}F_pA$ of $A$, both are intimately
related to the structure theory of $A$. Let $\prec_{gr}$ be an
$\NZ$-graded monomial ordering on the standard $K$-basis $\B$ of
$\KS$, and let $\langle\LM (I)\rangle$ be the monomial ideal of
$\KS$ generated by the set $\LM (I)$ of leading monomials of $I$
with respect to $\prec_{gr}$. If $\Omega$ is the unique reduced
monomial generating set of $\langle\LM (I)\rangle$ such that the
graph $\Gamma_{\rm C}(\Omega )$ of $n$-chains (in the sense of [1],
[11]) does not contain any $d$-chains, then, based on ([2], Theorem
4), it was proved in [7], without any extra assumption, that
$$\begin{array}{l}\hbox{gl.dim}A\le ~\hbox{gl.dim}G^{\NZ}(A)\le~\hbox{gl.dim}(\KS /\langle\LM (I)\rangle )\le d,\\
\hbox{gl.dim}\T A\le d+1,\end{array}$$
where gl.dim abbreviates the phrase ``global dimension". This note
aims to strengthen the above results, that is,  firstly we will show
further that if $\KS /\langle\LM (I)\rangle$ has the polynomial
growth of degree $m$,  then\par
(i) the following two equalities hold:
$$\hbox{gl.dim}G^{\NZ}(A)=~\hbox{gl.dim}(\KS /\langle\LM
(I)\rangle )=m,\quad \hbox{gl.dim}\T A = m+1;$$   \par
(ii) $I$ is generated by a finite Gr\"obner basis $\G$, and
consequently both $G^{\NZ}(A)$ and $\T A$ are defined by finite
Gr\"obner bases;\par
and secondly, we demonstrate, by examining interesting examples,
that bringing the problem of determining (i) above down-to-earth, if
we start with (ii), i.e.,  with a finite  Gr\"obner basis for $I$,
then an effective solution to (i) may be achieved. Moreover, the
last example given in section 3 will indicate in passing that under
the assumption of ([2], Theorem 6), the Hilbert series of a finitely
presented monomial algebra does not always have the form
$\prod^d_{i=1}(1-z^{e_i})^{-1}$ as asserted in loc. cit.} \v5
Throughout this paper we let $\KS$ denote the free $K$-algebra
$K\langle X_1,...,X_n\rangle$,  and let $\B$ be the standard
$K$-basis of $\KS$ consisting of words in the alphabet $X=\{
X_1,...,X_n\}$ (including the empty word which gives the identity
element 1). Unless otherwise stated, the $\NZ$-gradation of $\KS$
means any positive weight $\NZ$-gradation of $\KS$ by assigning to
each  $X_i$ a  positive degree $n_i$, $1\le i\le n$. Moreover,
ideals mean two-sided ideals, and if $M\subset\KS$, then we use
$\langle M\rangle$ to denote the ideal of $\KS$ generated by $M$.
For a general theory on Gr\"obner bases in $\KS$, we refer to
[9].\v5

\section*{1. Preliminaries}
In this section we recall several well-known algorithmic results
from [10], [2], and [5], that will be used in deriving the main
results of this note.\v5
Adopting the commonly used terminology in computational algebra, let
us call elements in $\B$ the monomials. Given a monomial ordering
$\prec$ on $\B$, as usual we write $\LM (f)$ for the leading
monomial of $f\in\KS$; and if $S$ is any subset of $\KS$, then we
write $\LM (S)$ for the set of leading monomials of $S$, i.e., $\LM
(S)=\{ \LM (f)~|~f\in S\}$. If $\G$ is a Gr\"obner basis in $\KS$
with respect to $\prec$, then it is well-known that we may always
assume that $\G$ is {\it LM}-{\it reduced}, that is, $g_1$,
$g_2\in\G$ and $g_1\ne g_2$ implies $\LM (g_1)\not |~\LM (g_2)$.
Consequently, if $\Omega$ is a subset of $\B$ satisfying  $u_i\not
|~u_j$ for all $u_i$, $u_j\in\Omega$ with $i\ne j$, then we just say
simply that $\Omega$ {\it is reduced}.\par
Let $\Omega =\{ u_1, ...,u_s\}$ be a reduced finite subset of $\B$.
For each $u_i\in\Omega$, say $u_i=X_{i_1}^{\alpha_1}\cdots
X_{i_s}^{\alpha_s}$ with $X_{i_j}\in X$ and $\alpha_j\in\NZ$,  we
write $l(u_i)=\alpha_1+\cdots +\alpha_s$ for the length of $u_i$ .
Put
$$\ell =\max\{ l(u_i)~|~u_i\in\Omega\} .$$
 Then the {\it Ufnarovski graph}  of $\Omega$ (introduced by
V. Ufnarovski in [10]), denoted $\Gamma (\Omega)$, is defined as a
directed graph, in which the set of vertices $V$ is given by
$$V=\{ v_i~|~v_i\in\B -\langle\Omega\rangle,~l(v_i)=\ell -1\} ,$$
and the set of edges $E$ contains the edge $v_i\rightarrow v_j$ if
and only if there exist $X_i$, $X_j\in X$ such that
$v_iX_i=X_jv_j\in\B -\langle\Omega\rangle$. Thus, for an LM-reduced
finite Gr\"obner basis $\G =\{ g_1,...,g_s\}$ in $\KS$, the
Ufnarovski graph of $\G$ is defined to be the Ufnarovski graph
$\Gamma (\LM (\G))$ of the reduced subset of monomials $\LM (\G )=\{
\LM (g_1),...,$ $\LM (g_s)\}$.{\parindent=0pt\v5
{\bf Remark} To better understand the practical application of
$\Gamma (\Omega )$, it is essential to notice that the Ufnarovski
graph is defined by using the {\it length $l(u)$ of the monomial}
({\it word}) $u\in \B$ instead of using the {\it degree of $u$ as an
$\NZ$-homogeneous element in $\KS$}, though both notions coincide
when each $X_i$ is assigned to degree 1.}\v5
The first effective application of $\Gamma (\Omega )$ was made to
determine the growth of the monomial algebra $\KS
/\langle\Omega\rangle$. {\parindent=0pt\v5
{\bf 1.1. Theorem} ([10], 1982) Let $\Omega=\{ u_1,...,u_s\}$ be a
reduced finite subset of $\B$, and let $\Gamma (\Omega )$ be the
Ufnarovski graph of $\Omega$ as defined above. Then the growth of
$\KS /\langle\Omega\rangle$ is alternative. It is exponential (i.e.,
the Gelfand-kirillov dimension of $\KS /\langle\Omega\rangle$ is
$\infty$) if and only if there are two different cycles  with a
common vertex in the graph $\Gamma (\Omega )$; Otherwise, $\KS
/\langle\Omega\rangle$ has the polynomial growth of degree $m$
(i.e.,  the Gelfand-Kirillov dimension of $\KS
/\langle\Omega\rangle$ is equal to $m$), where $m$ is, among all
routes of $\Gamma (\Omega )$,  the largest number of distinct cycles
occurring in a single route.\par\QED}\v5
Let the free $K$-algebra $\KS =K\langle X_1,...,X_n\rangle$  be
equipped with the augmentation map $\varepsilon$ sending each $X_i$
to zero. For any ideal $J$ contained in the augmentation ideal
$\langle X_1,...,X_n\rangle$ (i.e., the kernel of $\varepsilon$), by
using the $n$-chains determined by $\LM (J)$, D. J. Anick
constructed in [1] a free resolution of the trivial module $K$ over
the quotient algebra $\KS /J$ which gave rise to several efficient
applications to the homological aspects of associative algebras
([1], [2]). Let $\Omega\subset \B$ be a reduced (finite or infinite)
subset of monomials. Following [1] and [2], V. Ufnarovski
constructed in [11] the {\it graph of $n$-chains} of  $\Omega$ as a
directed graph $\Gamma_{\rm C}(\Omega )$, in which the set of
vertices $V$ is defined as
$$V=\{ 1\}\cup X\cup\{\hbox{all proper suffixes of}~u\in\Omega\},$$
and the set of edges $E$ consists of all edges
$$1~\mapright{}{}~X_i~\hbox{for every}~X_i\in X$$
and edges defined by the rule: for $u,v\in V-\{ 1\}$,
$$\begin{array}{rcl} ~u~\mapright{}{}~v~\hbox{in}~E&\Leftrightarrow &
\hbox{there is a {\it unique}}~w=X_{i_1}\cdots X_{i_{m-1}}X_{i_m}\in\Omega \\
&{~~~~}&\hbox{such that}~uv=\left\{\begin{array}{l} w,~\hbox{or}\\
sw~\hbox{with}~s\in\B,~sX_{i_1}\cdots X_{i_{m-1}}\in\B
-\langle\Omega\rangle\end{array}\right.\end{array}$$ For $n\ge -1$,
an $n$-{\it chain} of $\Omega$ is a monomial (word) $v=v_1\cdots
v_nv_{n+1}$ given by a route of length $n+1$ starting from $1$ in
$\Gamma_{\rm C}(\Omega )$:
$$1\r v_1\r\cdots\r v_n\r v_{n+1}$$
Writing $C_n$ for the set of all $n$-chains of $\Omega $, it is
clear that $C_{-1}=\{ 1\}$, $C_0=X$, and $C_1=\Omega$.\par
For an LM-reduced (finite or infinite) Gr\"obner basis $\G$ of
$\KS$, $\Gamma_{\rm C}(\LM (G) )$ is referred to the graph of
$n$-chains of $\G$. {\parindent=0pt\v5 {\bf Remark} As with the
Ufnarovski graph $\Gamma (\Omega )$ defined in section 2, to better
understand the practical application of the graph $\Gamma_{\rm
C}(\Omega )$ of $n$-chains determined by $\Omega$ in the subsequent
sections and the next chapter, it is essential to notice that {\it
an $n$-chain is defined by a route of length $n+1$ starting with 1,
as described above,  instead of by the degree of the
$\NZ$-homogeneous element $v=v_1\cdots v_nv_{n+1}$ read out of that
route}. \v5
{\bf 1.2. Theorem} ([2], Theorem 4) Let $\Omega\subset \B$ be a
reduced subset of monomials. Then gl.dim$(\KS /\langle\Omega\rangle
)\le m$ if and only if the graph $\Gamma_{\rm C}(\Omega )$ of
$n$-chains of $\Omega$ does not contain any $m$-chains.\par\QED \v5
\v5 {\bf 1.3. Theorem} ([2], Theorem 6) Let $\Omega$ be as in
Theorem 1.2. Suppose $\KS /\langle\Omega\rangle$ has finite global
dimension $m$. If $\KS /\langle\Omega\rangle$ does not contain a
free subalgebra of two generators, then the following statements
hold.\par
(i) $\KS /\langle\Omega\rangle$ is finitely presented, that is,
$\langle\Omega\rangle$ is finitely generated.\par
(ii) $\KS /\langle\Omega\rangle$ has the polynomial growth of degree
$m$.\par
(iii) The Hilbert series of $\KS /\langle\Omega\rangle$ is of the
form $H_{\KS
/\langle\Omega\rangle}(t)=\prod^m_{i=1}(1-t^{e_i})^{-1}$, where each
$e_i$ is a positive integer, $1\le i\le m$.\par\QED}\v5
By using the above theorem, the following result was derived by T.
Gateva-Ivanova in [5].{\parindent=0pt\v5
{\bf 1.4. Theorem} ([5], Theorem II) Let $J$ be an $\NZ$-graded
ideal of $\KS$ and $R=\KS /J$ the corresponding $\NZ$-graded algebra
defined by $J$. Suppose that the associated monomial algebra $\OV
R=\KS /\langle\LM (J)\rangle$ of $R$ has finite global dimension and
the polynomial growth of degree $m$, where $\LM (J)$ is taken with
respect to a fixed $\NZ$-graded monomial ordering $\prec_{gr}$ on
$\B$ (see the definition in the next section). Then the following
statements hold.\par
(i) gl.dim$R=$ gl.dim$\OV R=m$.\par
(ii) The ideal $J$ has a finite Gr\"obner basis.\par
(iii) The Hilbert series of $R$ is of the form
$H_R(t)=\prod^m_{i=1}(1-t^{e_i})^{-1}$, where each $e_i$ is a
positive integer, $1\le i\le n$.\par\QED}\v5

\section*{2. The Main Results}
In this section we show how to obtain the results of (i) -- (ii)
announced in the beginning of this note.\v5
Concerning the first equality
$\hbox{gl.dim}G^{\NZ}(A)=~\hbox{gl.dim}(\KS /\langle\LM (I)\rangle
)=m$, let us first recall some results from [8], [6], and [7]. Note
that with respect to the fixed weight $\NZ$-gradation of $\KS$,
every  element $f\in\KS$ has a unique expression $f=\sum_{i=1}^pF_p$
with $F_p\ne 0$, where each $F_i$ is an $\NZ$-homogeneous element of
degree $i$ in $\KS$, $1\le i\le p$. We call $F_p$ the $\NZ$-{\it
leading homogeneous element} of $f$, denoted $\LH (f)$, i.e., $\LH
(f)=F_p$. For a subset $S\subset\KS$, we then write $\LH (S)$ for
the set of $\NZ$-leading homogeneous elements of $S$, that is, $\LH
(S)=\{\LH (f)~|~f\in S\}$.{\parindent=0pt\v5
{\bf 2.1. Proposition} ([8], Proposition 2.2.1; [7], Theorem 1.1)
Let $I$ be an arbitrary ideal of $\KS$ and $A=\KS /I$. Considering
the $\NZ$-filtration $FA$ of $A$ induced by the weight $\NZ$-grading
filtration $F\KS$ of $\KS$, if $G^{\NZ}(A)$ is the associated
$\NZ$-graded algebra of $A$ determined by $FA$, then there is an
$\NZ$-graded $K$-algebra isomorphism
$$\KS /\langle\LH (I)\rangle\cong G^{\NZ}(A).$$\par\QED}\v5
Also recall that any well-ordering $\prec$ on $\B$ may be used to
define a new ordering $\prec_{gr}$: for $u$, $v\in\B$,
$$\begin{array}{rcl} u\prec_{gr} v&\Leftrightarrow& d(u)<d(v)\\
&~~&\hbox{or}~d(u)=d(v)~\hbox{and}~u\prec v,\end{array}$$
where $d(~)$ is referred to the degree function on homogeneous
elements of $\KS$. If $\prec_{gr}$ is a monomial ordering on $\B$,
then it is clled an $\NZ$-{\it graded monomial ordering}. Typical
$\NZ$-graded monomial ordering on $\B$ is the well-known {\it
$\NZ$-graded (reverse) lexicographic ordering}. {\parindent=0pt\v5
{\bf 2.2. Theorem} ([8], Theorem 2.3.2; [7], Proposition 3.2.) Let
$I$ be an arbitrary ideal of $\KS$ and $\G\subset I$. Then $\G$ is a
Gr\"obner basis of $I$ with respect to some $\NZ$-graded monomial
ordering $\prec_{gr}$ on $\B$  if and only if the set of
$\NZ$-leading homogeneous elements $\LH (\G )$ of $\G$ is a
Gr\"obner basis for the $\NZ$-graded ideal $\langle\LH (I)\rangle$
with respect to $\prec_{gr}$.\par\QED}\v5
We are ready to obtain the first result of this
section.{\parindent=0pt\v5
{\bf 2.3. Theorem} Let $I$ be an arbitrary ideal of $\KS$, $A=\KS
/I$, and $\OV A=\KS /\langle\LM (I)\rangle$ the associated monomial
algebra of $A$ with respect to a fixed $\NZ$-graded monomial
ordering $\prec_{gr}$ on $\B$. Suppose that gl.dim$\OV A<\infty$,
and that $\OV A$ has the polynomial growth of degree $m$. With
notation as before, the following statements hold.\par
(i) gl.dim$G^{\NZ}(A)=$ gl.dim$\OV A=m$.\par
(ii) The ideal $I$ has a finite Gr\"obner basis $\G$, and $\LH (\G
)$ is a finite homogeneous Gr\"obner basis in $\KS$ such that
$G^{\NZ}(A)\cong\KS /\langle\LH (\G )\rangle$.\vskip 6pt
{\bf Proof} (i) Since we are using the $\NZ$-graded monomial
ordering $\prec_{gr}$, it is straightforward that $\LM (I)=\LM
(\langle\LH (I)\rangle )$. Hence, both algebras $A=\KS /I$ and $\KS
/\langle\LH (I)\rangle$ have the same associated monomial algebra
$\OV A=\KS /\langle\LM (I)\rangle$. So, by making use of the unique
reduced monomial generating set $\Omega$ of the monomial ideal
$\langle\LM (I)\rangle$, the equality in (i) follows from Theorem
1.3(ii), Theorem 1.4(i), and Proposition 2.1.\par
(ii) It is a well-known fact that if $\langle\LM (I)\rangle$ is
finitely generated, then $I$ is generated by a finite Gr\"obner
basis. That $I$ has a finite Gr\"obner basis $\G$ is guaranteed by
Theorem 1.3(i). Thus, the second assertion of (ii) concerning $\LH
(\G )$ follows from Proposition 2.1 and Theorem 2.2.\par\QED}\v5
It remains to show that under the same assumption as in Theorem 2.3,
the Rees algebra $\T A$ of $A$ has the properties listed in the
beginning of this paper. To this end, we need a little more
preparation.\par
Let $\prec_{gr}$ be a fixed $\NZ$-graded monomial ordering on $\B$
with respect to the positive weight $\NZ$-gradation of $\KS$. If
$f\in \KS$ has the linear presentation by elements of $\B$:
$$f=\LC (f)\LM (f)+\sum_i\lambda_iw_i~\hbox{with}~
\LM (f)\in\B\cap\KS_p,~w_i\in\B\cap\KS_{q_i} ,$$
where $\LC (f)$ is the leading coefficient of $f$, then $f$
corresponds to a unique homogeneous element in the free algebra
$K\langle X,T\rangle =K\langle X_1,...,X_n,T\rangle$, i.e., the
element
$$\T f=\LC (f)\LM (f)+\sum_i\lambda_iT^{p-q_i}w_i.$$
Assigning to $T$ the degree 1 in $K\langle X,T\rangle$ and using the
fixed positive weight of $X$ in $\KS$, we get the weight
$\NZ$-gradation of $K\langle X,T\rangle$ which extends the weight
$\NZ$-gradation of $\KS$. Consequently, writing $\T{\B}$ for the
standard $K$-basis of $K\langle X,T\rangle$, we may extend
$\prec_{gr}$ to an $\NZ$-graded monomial ordering
$\prec_{_{T\hbox{-}gr}}$ on $\T{\B}$ such that
$$T\prec_{_{T\hbox{-}gr}}X_i,~1\le i\le n,~\hbox{and hence}~
\LM (f)=\LM (\T f).$$ If $I$ is an ideal of $\KS$ generated by the
subset $S$, then we put
$$\begin{array}{l} \T I=\{ \T f~ |~f\in I\}\cup\{ X_iT-TX_i~ |~1\le i\le n\} ,\\
\\
\T{S}=\{ \T f~ |~f\in S\}\cup\{ X_iT-TX_i~ |~1\le i\le n\}
.\end{array}$$
The assertion of the next proposition was inferred in ([8], Theorem
2.3.1; [6], CH.III, Corollary 3.8) and a detailed proof was given in
([7], section 8){\parindent=0pt\v5
{\bf 2.4. Proposition} With the preparation made above, let $\G$ be
a Gr\"obner basis of the ideal $I$ in $\KS$ with respect to
$\prec_{gr}$, and $A=\KS /I$. Then $\T{\G}$ is a  Gr\"obner basis
for the $\NZ$-graded ideal $\langle\T I\rangle$ in $K\langle
X,T\rangle$ with respect to  $\prec_{_{T\hbox{-}gr}}$, and
consequently  $\T A\cong K\langle X,T\rangle /\langle\T I\rangle
=K\langle X,T\rangle/\langle\T{\G}\rangle$, where $\T A$ is the Rees
algebra of $A$ defined by the $\NZ$-filtration $FA$ induced by the
weight $\NZ$-grading filtration $F\KS$ of $\KS$.\par\QED}\v5
Now, the following result is established.{\parindent=0pt\v5
{\bf 2.5. Theorem} Let $I$ be an arbitrary ideal of $\KS$, $A=\KS
/I$, and $\OV A=\KS /\langle\LM (I)\rangle$ the associated monomial
algebra of $A$ with respect to a fixed $\NZ$-graded monomial
ordering $\prec_{gr}$ on $\B$. Suppose that gl.dim$\OV A<\infty$,
and that $\OV A$ has the polynomial growth of degree $m$. With
notation as before, the following statements hold.\par
(i) gl.dim$\T A=m+1$.\par
(ii) The ideal $I$ has a finite Gr\"obner basis $\G$, and $\T{\G }$
is a finite homogeneous Gr\"obner basis of $K\langle X,T\rangle$
such that  $\T A\cong K\langle X,T\rangle
/\langle\T{\G}\rangle$.\vskip 6pt
{\bf Proof} By the assumption, $\OV A$ has Gelfand-Kirillov
dimension $m$. Hence $A$ has Gelfand-Kirillov dimension $m$. It
follows from ([7], Theorem 8.3(i)) that $\T A$ has Gelfand-Kirillov
dimension $m+1$. But by the assumption and Theorem 2.3, $I$ has a
finite Gr\"obner basis $\G$. So by Proposition 2.4,  $\T{\G}$ is a
finite Gr\"obner basis for the ideal $\langle\T I\rangle$ and $\T
A\cong K\langle X,T\rangle /\langle\T{\G}\rangle$. Thus, $\T A$ must
have the polynomial growth of degree $m+1$ by Theorem 1.1.
Therefore, the associated monomial algebra $K\langle X,T\rangle
/\langle\LM (\T{\G})\rangle$ of $\T A$ has the polynomial growth of
degree $m+1$. In order to finish the proof, by Theorem 1.4, it
remains to show that $K\langle X,T\rangle /\langle\LM
(\T{\G})\rangle$ has finite global dimension. To see this clearly,
we quote the argument from the proof of ([7], Theorem 8.5) as
follows.}\par
Note that by the definition of $\prec_{_{T\hbox{-}gr}}$ we have
$$\LM (\widetilde{\G})=\{ \LM (g),~X_iT~|~g\in\G,~1\le i\le n\} .$$
Thus, it is easy to see that the graph of $n$-chains $\Gamma_{\rm
C}(\LM (\G ))$ of $\G$ is a subgraph of the graph of $n$-chains
$\Gamma_{\rm C}(\LM (\widetilde{\G}))$ of $\widetilde{\G}$, and that
{\parindent=.7truecm
\item{(a)} the graph $\Gamma_{\rm C}(\LM
(\widetilde{\G}))$ of $\widetilde{\G}$ has no edge of the form $T\r
v$ for all $v\in \widetilde{V}$, where $\widetilde{V}$ is the set of
vertices of $\Gamma_{\rm C}(\LM (\widetilde{\G}))$;
\item{(b)} if $v\in\widetilde{V}$ is of the form $v=sX_j$, $s\in\B$, then $\Gamma_{\rm C}(\LM
(\widetilde{\G}))$ contains the edge $v\r T$;
\item{(c)} any $d+1$-chain in $\Gamma_{\rm
C}(\LM (\widetilde{\G}))$ is of the form
$$1\r X_i\r v_1\r v_2\r\cdots\r v_{d-1}\r T,$$
where
$$1\r X_i\r v_1\r v_2\r\cdots\r v_{d-1}$$
is a $d$-chain in $\Gamma_{\rm C}(\LM (\G ))$.\parindent=0pt\par
Therefore,  if the graph $\Gamma_{\rm C}(\LM (\G ))$ does not
contain any $d$-chain, then $\Gamma_{\rm C}(\LM (\T{\G}))$ does not
contain any $d+1$-chain. Hence, if gl.dim$\OV A=\KS /\langle\LM
(I)\rangle <\infty$, then gl.dim$K\langle X,T\rangle /\langle\LM
(\T{\G})\rangle <\infty$ by Theorem 1.2, as desired.\QED}\v5

\section*{3. Examples of Calculating gl.dim$G^{\NZ}(A)$ and gl.dim$\T A$}
Let $I$, $A=\KS /I$, $\OV A=\KS /\langle\LM (I)\rangle$,
$G^{\NZ}(A)$, and $\T A$ be as in section 2. Combining Theorem 1.1,
Theorem 1.2, Proposition 2.1 and Proposition 2.4, it is now clear
that if, with respect to a fixed $\NZ$-graded monomial ordering
$\prec_{gr}$ on $\B$, we start with a finite Gr\"obner basis $\G =\{
g_1,...,g_s\}$ for the ideal $I$, then the equalities of Theorem
2.3(i) and Theorem 2.5(i) may be determined in a computational
way:{\parindent=0pt\vskip 6pt
(1) Determine whether $\OV A$ has polynomial growth by checking the
Ufnarovski graph $\Gamma (\LM (\G ))$; if $\OV A$ has polynomial
growth, then the degree $m$ is read out of the graph simultaneously.
\par
(2) Determine whether $\OV A$ has finite global dimension by
checking the set $C_n$ of $n$-chains in the graph $\Gamma_{\rm
C}(\LM (\G ))$; if $C_d=\emptyset$ for some $d$, then gl.dim$\OV
A\le d$.\par
(3) If both (1) and (2) have a positive result, i.e., $\OV A$ has
the polynomial growth of degree $m$ and gl.dim$\OV A<\infty$, then
immediately we can write down the following:
$$\hbox{gl.dim}G^{\NZ}(A)=m,\quad\hbox{gl.dim}\T A=m+1.$$}\par
Let us point out incidentally that when the above (2) is done, the
Hilbert series for both $G^{\NZ}(A)$ and $\T A$ may also be written
down just by using the $n$-chains in $\Gamma_{\rm C}(\LM (\G ))$,
respectively the $n$-chains in $\Gamma_{\rm C}(\LM (\T{\G}))$ as
described in the proof of Theorem 2.5, and ([1], formula 16), that
is,
$$H_{G^{\NZ}(A)}(t)=\left (1-\sum_{i=0}(-1)^iH_{C_i}(t)\right )^{-1},\quad
H_{\T A}(t)=\left (1-\sum_{i=0}(-1)^iH_{\T C_i}(t)\right )^{-1},$$
where $H_{C_i}(t)$ denotes the Hilbert series of the $\NZ$-graded
$K$-module spanned by the set $C_i$ of $i$-chains in $\Gamma_{\rm
C}(\LM (\G ))$, and similarly, $H_{\T C_i}(t)$ denotes the Hilbert
series of the $\NZ$-graded $K$-module spanned by the set $\T C_i$ of
$i$-chains in $\Gamma_{\rm C}(\LM (\T{\G} ))$.\vskip 6pt
We illustrate the computational procedure mentioned above by
examining several examples. Notations are maintained as before. \v5
Let $\Omega =\{ X_jX_i~|~1\le i<j\le n\}\subset\B\subset\KS$. Then
by ([2], Example 3), gl.dim$(\KS /\langle\Omega\rangle =n=$ the
degree of the polynomial growth of $\KS /\langle\Omega\rangle$. Note
that $\B -\langle\Omega\rangle =\{
X_1^{\alpha_1}X_2^{\alpha_2}\cdots
X_n^{\alpha_n}~|~\alpha_1,...,\alpha_n\in\NZ\}$, which gives rise to
a PBW $K$-basis for $\KS /\langle\Omega\rangle$. Enlightened by this
fact, our first example will be the algebra $A=\KS /I$, which, with
respect to a fixed monomial ordering $\prec$ on $\B$, has the
property that $\B -\LM (I) =\{ X_1^{\alpha_1}X_2^{\alpha_2}\cdots
X_n^{\alpha_n}~|~\alpha_1,...,\alpha_n\in\NZ\}$ yields a PBW
$K$-basis for both $A$ and $\KS /\langle\LM (I)\rangle$. Below let
us describe first the ideal $I$ by Gr\"obner basis (probably a known
result but the author has lack of a proper
reference).{\parindent=0pt\v5
{\bf 3.1. Proposition} Let $I$ be an ideal of $\KS$ and $A=\KS /I$.
The following two statements are equivalent with respect to a fixed
monomial ordering $\prec$ on $\B$:\par
(i) $\B -\LM (I)= \{ X_1^{\alpha_1}X_2^{\alpha_2}\cdots
X_n^{\alpha_n}~|~\alpha_j\in\NZ\}$ and hence $A$ has the PBW
$K$-basis $\{ \OV{X_1}^{\alpha_1}\OV{X_2}^{\alpha_2}\cdots
\OV{X_n}^{\alpha_n}~\Big |~\alpha_j\in\NZ\}$, where each $\OV{X_i}$
is the image of $X_i$ in $A$;\par
(ii) $I$ is generated by a reduced Gr\"obner basis of the form
$$\begin{array}{rcl} \G &=&\left\{ R_{ji}=X_jX_i-F_{ji}~\Big |~F_{ji}\in\KS ,~1\le i<j\le n\right\}\\
&{~}&\hbox{satisfying}~\LM (R_{ji})=X_jX_i,~1\le i<j\le
n,~\hbox{and}\\
&{~}&F_{ji}=\SUM^m_{p=1}\lambda_pw_p~\hbox{with}~\lambda_p\in
K^*,~w_p\in\{ X_1^{\alpha_1}X_2^{\alpha_2}\cdots
X_n^{\alpha_n}~|~\alpha_j\in\NZ\} .\end{array}$$  \vskip 6pt
{\bf Proof} (i) $\Rightarrow$ (ii) First, we show that under the
assumption of (i) $I$ has a finite Gr\"obner basis $\G=\{
R_{ji}~|~1\le i<j\le n\}$ of the described form such that $\LM
(R_{ji})=X_jX_i$. By classical G\"obner basis theory, it is
sufficient to prove that the reduced monomial generating set
$\Omega$ of $\langle\LM (I)\rangle$ is of the form $\Omega =\{
X_jX_i~|~1\le i<j\le n\} .$ To see this, recall  that
$$\Omega =\{ w\in\LM (I)~~|~\hbox{if}~u\in\LM (I)~\hbox{and}~u|w~\hbox{then}~u=w\} ,$$
and consequently $\B -\LM (I)$ is obtained by the division by
$\Omega$. Since $X_i\in\B -\LM (I)$,  and $X_jX_i\not\in\B -\LM
(I)$, it follows that $X_jX_i\in\Omega$, $1\le i<j\le n$.
Furthermore, noticing the feature of monomials in $\B -\LM (I)$, it
is clear that the only monomials of length 2 contained in $\Omega$
are $X_jX_i$, $1\le i<j\le n$, and that $\Omega$ cannot contain
monomials of length $\ge 3$. Therefore, $\Omega$ has the form as we
claimed above, and $\Omega$ determines a Gr\"obner basis $\G =\{
R_{ji}~|~1\le i<j\le n\}$ with $\LM (R_{ji})=X_jX_i$. It follows
from classical Gr\"obner basis theory that $\G$ can be reduced
further to a Gr\"obner basis such that $F_{ji}=R_{ji}-\LM (R_{ji})$
is a normal element, which is clearly a linear combination of
monomials in $\{ X_1^{\alpha_1}X_2^{\alpha_2}\cdots
X_n^{\alpha_n}~|~\alpha_j\in\NZ\}$, but this means actually that
$\G$ is a reduced Gr\"obner bsis.{\parindent=0pt\par (ii)
$\Rightarrow$ (i) If $\G$ is a Gr\"obner basis of $I$ as described,
then since $\langle\LM (I)\rangle =\langle\LM (\G )\rangle$ and $\LM
(R_{ji})=X_jX_i$, $1\le i<j\le n$, the division by $\LM (\G)$ yields
$$\B -\LM (I)=\left\{ X_1^{\alpha_1}X_2^{\alpha_2}\cdots X_n^{\alpha_n}~\Big |~\alpha_i\in\NZ\right\} ,$$
as desired.\QED}\v5
{\bf Remark} The last proposition tells us that in order to have a
PBW $K$-basis in terms of Gr\"obner basis in $\KS$, it is necessary
to consider a finite subset $\G$ of $\KS$ as described in
Proposition 3.1(ii).}\v5
In light of Proposition 3.1, the next result is now obtained by
using Theorem 2.3, Theorem 2.5 and ([2], Example 3).
{\parindent=0pt\v5
{\bf 3.2. Theorem} Fixing a positive weight $\NZ$-gradation for the
free $K$-algebra $\KS =K\langle X_1,...,X_n\rangle$, let
$\prec_{gr}$ be an $\NZ$-graded monomial ordering on the standard
$K$-basis $\B$ of $\KS$. With notation as in section 2, if $\G$ is a
Gr\"obner basis in $\KS$ with respect to $\prec_{gr}$, such that the
ideal $I=\langle\G\rangle$ has the property that $\B -\LM (I)=\{
X_1^{\alpha_1}X_2^{\alpha_2}\cdots
X_n^{\alpha_n}~|~\alpha_1,...,\alpha_n\in\NZ\}$, then with respect
to the $\NZ$-filtration $FA$ of $A=\KS /I$ induced by the weight
$\NZ$-grading filtration $F\KS$,
$$\begin{array}{l} \hbox{gl.dim}G^{\NZ}(A)=~\hbox{gl.dim}(\KS /\langle\LH (\G )\rangle) =
\hbox{gl.dim}(\KS /\langle\LM (\G )\rangle) =n;\\
\hbox{gl.dim}\T A=~\hbox{gl.dim}(K\langle X,T\rangle /\langle\LM
(\T{\G})\rangle) =n+1.\end{array}$$\par\QED \v5
{\bf Example} (1) Rather than quoting those well-known Gr\"obner
bases that give rise to a PBW $K$-basis, let us look at a small one.
Consider the ideal $I=\langle R_{21}\rangle$ of the free $K$-algebra
$\KS =K\langle X_1,X_2\rangle$ generated by the single element
$$R_{21}=X_2X_1-qX_1X_2-\alpha X_2-f(X_1),$$
where $q$, $\alpha\in K$, and $f(X_1)$ is a polynomial in the
variable $X_1$. Assigning to $X_1$ the degree 1, then in either of
the following two cases:\par (a) deg$f(X_1)\le 2$, and $X_2$ is
assigned to the degree 1;\par       (b) deg$f(X_1)=n\ge 3$, and
$X_2$ is assigned to the degree n,\par $\G =\{ R_{21}\}$ forms a
Gr\"obner basis for $I$ with respect to the $\NZ$-graded
lexicographic ordering $X_1\prec_{gr} X_2$, such that $\LM (\G )=\{
X_2X_1\}$. Putting $A=K\langle X_1,X_2\rangle /I$, and noticing that
in both gradations $\LH (R_{21})= X_2X_1-qX_1X_2$ and
$\T{R_{21}}=X_2X_1-qX_1X_2-\alpha TX_2-\T{f(X_1)}$, it follows from
Theorem 3.2 that gl.dim$G^{\NZ}(A)=2$, gl.dim$\T A=3$. \v5
{\bf Example} (2) Let the free $K$-algebra $\KS =K\langle
X_1,X_2\rangle$ be equipped with a  positive weight $\NZ$-gradation,
such that $d(X_1)=n_1$ and $d(X_2)=n_2$,  and let $\G=\{ g_1,g_2\}$
be any  Gr\"obner basis with respect to some $\NZ$-graded monomial
ordering $\prec_{gr}$ on the standard basis $\B$ of $\KS$, such that
$\LM (g_1)=X_1^2X_2$ and $\LM (g_2)=X_1X_2^2$. For instance, with
respect to the $\NZ$-graded lexicographic ordering $X_2\prec_{gr}
X_1$,
$$\begin{array}{l} g_1=X_1^2X_2-\alpha X_1X_2X_1-\beta X_2X_1^2-\lambda X_2X_1-\gamma X_1,\\
g_2=X_1X_2^2-\alpha X_2X_1X_2-\beta X_2^2X_1-\lambda X_2^2-\gamma
X_2.\end{array} \quad\alpha ,\beta ,\lambda ,\gamma\in K.$$
Consider the algebra $A=\KS /\langle\G\rangle$ which is equipped
with the $\NZ$-filtration $FA$ induced by the weight $\NZ$-grading
filtration $F\KS$. With notation as before, the following statements
hold.\par (i) All three algebras $A$, $G^{\NZ}(A)$ and $\KS
/\langle\LM (\G )\rangle$ have the polynomial growth of degree 3,
while $\T A$ and $K\langle X_1,X_2,T\rangle /\langle\LM
(\T{\G})\rangle$ have the polynomial growth of degree 4.\par
(ii) gl.dim$G^{\NZ}(A)=$ gl.dim$\KS /\langle\LM (\G )\rangle =3$,
and gl.dim$\T A =$ gl.dim$K\langle X_1,X_2,T\rangle /\langle\LM
(\T{\G})\rangle = 4$.\par
In particular,  all results mentioned above hold for the down-up
algebra $A(\alpha ,\beta ,\gamma )$ (in the sense of [3]) which is
defined by the relations
$$\begin{array}{l} g_1=X_1^2X_2-\alpha X_1X_2X_1-\beta X_2X_1^2-\gamma X_1,\\
g_2=X_1X_2^2-\alpha X_2X_1X_2-\beta X_2^2X_1-\gamma X_2.\end{array}
\quad\alpha ,\beta ,\gamma\in K.$$ {\parindent=0pt\par
{\bf Proof} Put $\Omega =\{ X_1^2X_2,~X_1X_2^2\}$. Then the
Ufnarovski graph $\Gamma (\Omega )$ of $\Omega$ is presented by
$$
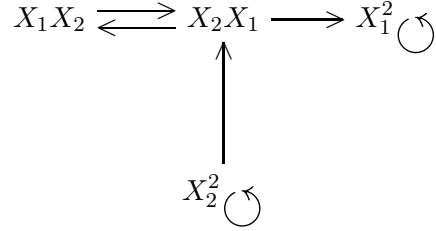
\begin{diagram} X_1X_2&\pile{\rTo \\ \lTo}&X_2X_1&\rTo &X_1^2{_{\hbox{\huge$\circlearrowleft$}}}\\
 &&\uTo&&\\
&&X_2^2{_{\hbox{\huge$\circlearrowleft$}}}&&\end{diagram}$$
which shows that $\KS /\langle\LM (\G )\rangle$ has the polynomial
growth of degree 3. Hence, by referring to the proof of Theorem 2.5,
the assertions of (i) are determined. Next, the graph of $n$-chains
$\Gamma_{\rm C}(\Omega )$ of $\Omega$ is presented by
$$\begin{array}{ccccl} &&&&X_2^2\\
&&&\nearrow&\\
&&X_1&&\\
&\nearrow&&\searrow&\\
1&&&&X_1X_2\\
&\searrow&&\swarrow&\\
&&X_2&&\end{array}$$ which shows that
$$C_{i-1}=\left\{\begin{array}{ll}
\{ X_1,X_2\},&i=1,\\
\{ X_1X_2^2,X_1^2X_2\},&i=2,\\
\{X_1^2X_2^2\},&i=3,\\ \emptyset ,&i\ge 4.\end{array}\right.$$
Note that $\LM (\T{\G})=\{X_1^2X_2,~X_1X_2^2,~X_1T,~X_2T\}$. By
referring to the proof of Theorem 2.5, it is then straightforward
that the set $\T C_{i-1}$ consisting of $i-1$-chains from
$\Gamma_{\rm C}(\LM (\T\G ))$ is
$$\T C_{i-1}=\left\{\begin{array}{ll}
\{ X_1,X_2,T\}&i=1,\\ \{ X_1T,X_2T,X_1X_2^2,X_1^2X_2\},&i=2,\\
\{ X_1^2X_2^2,X_1X_2^2T,X_1^2X_2T\},&i=3,\\
\{ X_1^2X_2^2T\},&i=4,\\ \emptyset ,&i\ge 5.\end{array}\right.$$
So the conditions of Theorem 2.3 and Theorem 2.5 are satisfied, and
consequently the assertions of (ii) are
determined.\QED{\parindent=0pt\v5 \def\HL{{\rm LH}}
{\bf Example} (3) Let the free $K$-algebra $\KS =K\langle
X_1,X_2\rangle$ be equipped with a positive weight $\NZ$-gradation
such that $d(X_1)=n_1$ and $d(X_2)=n_2$, and for any positive
integer $n$, let $\G =\{ g\}$ with $g=X_2^nX_1-qX_1X_2^n-F$ with
$q\in K$, $F\in\KS$. Consider the algebra $A=\KS /\langle\G\rangle$
which is equipped with the $\NZ$-filtration $FA$ induced by the
weight $\NZ$-grading filtration $F\KS$. With notation as before, the
following statements hold.\par
(i) All three algebras $A$, $G^{\NZ}(A)$ and $\KS /\langle\LM (\G
)\rangle$ have the polynomial growth of degree 2, while $\T A$ and
$K\langle X_1,X_2,T\rangle /\langle\LM (\T\G )\rangle$ have the
polynomial growth of degree 3.\par
(ii) gl.dim$G^{\NZ}(A)=$ gl.dim$\KS /\langle\LM (\G )\rangle =2$,
and gl.dim$\T A =$ gl.dim$K\langle X_1,X_2,T\rangle /\langle\LM
(\T\G  )\rangle = 3$.\par
(iii) If $n_1=n_2=1$, then $G^{\NZ}(A)$, respectively $\T A$, has
Hilbert series
$$H_{G^{\NZ}(A)}(t)=\frac{1}{1-2t+t^{n+1}}, ~\hbox{respectively}~H_{\T
A}(t)=\frac{1}{1-3t+2t^2+t^{n+1}-t^{n+2}},$$
which, in the case of $n\ge 2$, cannot be always the form
$$\frac{1}{(1-t^{e_1})(1-t^{e_2})},~\hbox{respectively}~\frac{1}{(1-t^{e_1'})(1-t^{e_2'})(1-t^{e_3'})}$$
as claimed in ([2], Theorem 6).\vskip 6pt
{\bf Proof} Consider the $\NZ$-graded lexicographic ordering
$X_1\prec_{gr} X_2$ with respect to the natural $\NZ$-gradation of
$\KS$. For any positive integer $n$, if $g=X_2^nX_1-qX_1X_2^n-F$
with $q\in K$ and $F\in\KS$ such that $\LM (F)\prec_{gr} X_2^nX_1$,
then it is easy to check that $\G =\{ g\}$ forms a Gr\"obner basis
for the ideal $I=\langle\G\rangle$.  Put $\Omega =\{\LM (\G
)=X_2^nX_1\}$. Then it is straightforward to verify that for $n=1$,
the Ufnarovski graph $\Gamma (\Omega )$ is presented by
$$
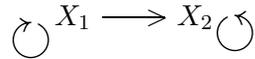
\begin{diagram} {~}_{_{\hbox{\huge$\circlearrowright$}}}X_1\rTo
X_2{_{\hbox{\huge$\circlearrowleft$}}}\end{diagram}$$
for $n=2$, the Ufnarovski graph $\Gamma (\Omega )$ is presented by
$$
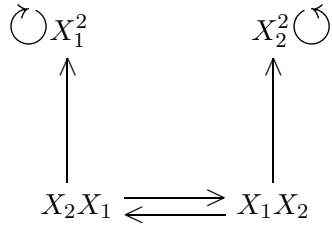
\begin{diagram} {{\hbox{\huge$\circlearrowright$}}}X_1^2~~~~&&~~~~X_2^2{{\hbox{\huge$\circlearrowleft$}}}\\
{~~}\uTo&&\uTo\\
{~~}X_2X_1&\pile{\rTo\\ \lTo}&X_1X_2\end{diagram}$$
and for $n\ge 3$, the Ufnarovski graph $\Gamma (\Omega )$ is
presented by \def\l{\leftarrow}
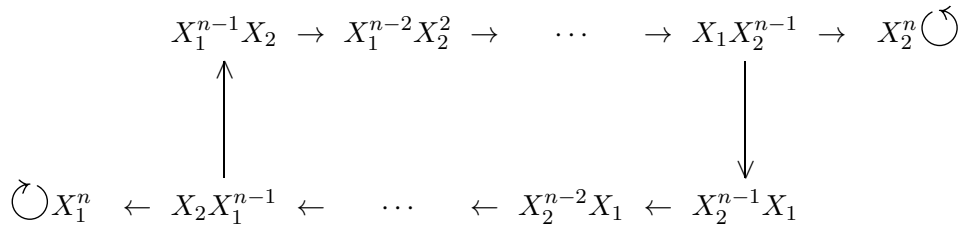
$$\begin{diagram}
&&X_1^{n-1}X_2&\r&X_1^{n-2}X_2^2&\r&\cdots&\r&X_1X_2^{n-1}&\r&X_2^n\hbox{\huge$\circlearrowleft$}\\
&&\uTo&&&&&&\dTo&&\\
\hbox{\huge$\circlearrowright$}X_1^n&\l&X_2X_1^{n-1}&\l&\cdots&\l&X_2^{n-2}X_1&\l&X_2^{n-1}X_1&&&\end{diagram}$$
Hence $\KS /\langle\LM (\G )\rangle$ has the polynomial growth of
degree 2, and then (i) follows. Since the graph $\Gamma_{\rm
C}(\Omega )$ of $n$-chains of $\Omega$ is presented by
$$\begin{array}{ccc} {~}&&{~}\\
{~}&&{~}\\
X_1&\longleftarrow&X_2\\ &\nwarrow~~~~~~~~\nearrow&\\
&1&\\
&(n=1)&\end{array}\quad\quad\quad \begin{array}{cccccccc} &&X_2^qX_1&&&&X_2^{n-1}X_1\\
&&&&&\nearrow&\\
X_1&&&&X_2&&\\
&\nwarrow&&\nearrow&&&\\
&&1&&&&\\
&&(n\ge 2)&&&&\end{array}$$
where $q\le n-2$ for the vertices $X_2^qX_1$,  it is clear that
$$C_{i-1}=\left\{\begin{array}{ll} \{ X_1,~X_2\},&i=1,\\ \{ X_2^nX_1\},&i=2,\\
~\emptyset ,&i\ge 3.\end{array}\right.$$
Also note that $\LM (\T\G )=\{ X_2^nX_1, X_2T, X_1T\}$. By referring
to the proof of Theorem 2.5, it is then straightforward that the set
$\T C_{i-1}$ consisting of $i-1$-chains from $\Gamma_{\rm C}(\LM
(\T\G ))$ is
$$\T C_{i-1}=\left\{\begin{array}{ll}
\{ X_1,X_2,T\}&i=1,\\ \{ X_1T,X_2T,X_2^nX_1\},&i=2,\\
\{ X_2^nX_1T\},&i=3,\\
\emptyset ,&i\ge 4.\end{array}\right.$$
So the conditions of Theorem 2.3 and Theorem 2.5 are satisfied, and
consequently the assertions of (ii) and (iii) are determined.\QED
\v5
{\bf Remark} For every positive integer $N=n+1\ge 2$, it was shown
in [7] that if $q\ne 0$ and $F\ne 0$ with total degree $<n$, then
the algebras defined by $\G$ and $\LH (\G )$ in the last example
provide a (non-monomial) homogeneous $N$-Koszul algebra and a
non-homogeneous $N$-Koszul algebra (in the sense of [4]),
respectively. \v5
\centerline{References}
\parindent=.8truecm

\item{[1]} D. J. Anick, On the homology of associative algebras, {\it
Trans. Amer. Math. Soc}., 2(296)(1986), 641--659.
\item{[2]} D. J. Anick, On Monomial algebras of finite global
dimension, {\it Trans. Amer. Math. Soc}., 1(291)(1985), 291--310.
\item{[3]} G.~Benkart, Down-up algebras and Witten's deformations of
the universal enveloping algebra of $sl_2$, {\it Contemp. Math.},
224(1999), 29--45.
\item{[4]} R. Berger and V. Ginzburg, Higher symplectic reflection
algebras and nonhomogeneous $N$-Koszul property, {\it J. Alg.},
1(304)(2006), 577--601.
\parindent=.8truecm
\item{[5]} T.~Gateva-Ivanova, Global dimension of associative
algebras, in: {\it Proc. AAECC-6}, LNCS, Vol.357, Springer-Verlag,
1989, 213--229.
\item{[6]} H. Li, {\it Noncommutative Gr\"obner Bases and
Filtered-Graded Transfer}, LNM, 1795, Springer-Verlag, 2002.
\item{[7]} H. Li, $\Gamma$-leading homogeneous algebras and Gr\"obner bases, {\it Alg. Colloq.},
to appear. (see a full version at
http://www.geocities.com/huishipp/gam$_{-}$grob.pdf )
\item{[8]} H. Li, Y. Wu and J. Zhang, Two applications of
noncommutative Gr\"obner bases, {\it Ann. Univ. Ferrara - Sez. VII -
Sc. Mat.}, XLV(1999), 1--24. (a full version is at
http://www.springerlink.com/content/v67351x2031p1079/)
\item{[9]} T. Mora, An introduction to commutative and noncommutative
Gr\"obner Bases, {\it Theoretic Computer Science}, 134(1994),
131--173.
\item{[10]} V. Ufnarovski, A growth criterion for graphs and algebras
defined by words, {\it Mat. Zametki}, 31(1982), 465--472 (in
Russian); English translation: {\it Math. Notes}, 37(1982),
238--241.
\item{[11]} V. Ufnarovski, On the use of graphs for computing a basis,
growth and Hilbert series of associative algebras, (in Russian
1989), {\it Math. USSR Sbornik}, 11(180)(1989), 417-428.

\end{document}